\documentclass[twoside,11pt]{article}
\usepackage{graphicx}
\usepackage{amssymb}
\usepackage{amsthm}
\newfont{\bbb}{msbm10 scaled\magstep 1}
\newcommand{\bd}{{\rm bd}}
\newcommand{\inter}{{\rm int}}

\font\bigbold=cmbx10 at 14 pt
\font\bbigbold=cmbx10 at 17.2 pt 
\date{}

\title
{\bbigbold The centroid Banach-Mazur distance between the parallelogram and the triangle}

\begin{document}

\baselineskip 17.5pt 

\maketitle

\vskip -1.3cm 
\centerline
{\bigbold Marek Lassak}

\vskip0.25cm
\pagestyle{myheadings} \markboth{Marek Lassak}{Centroid Banach-Mazur distance between the parallelogram and the triangle}

\noindent
{\bf Abstract.}
Let $C$ and $D$ be convex bodies in the Euclidean space $E^d$.
We define the centroid Banach-Mazur distance $\delta_{BM}^{\rm cen} (C, D)$ similarly to the classic Banach-Mazur distance $\delta_{BM} (C, D)$, but with the extra requirement that the centroids of $C$ and an affine image of $D$ coincide. 
We prove that for the parallelogram $P$ and the triangle $T$ in $E^2$ we have $\delta_{BM}^{\rm cen} (P, T) = \frac{5}{2}$.



\date{}

\maketitle

\section{Introduction}

The classical definition of the Banach-Mazur distance of centrally symmetric convex bodies of the Euclidean $d$-space $E^d$ is given by Banach \cite{[B]} in behalf of him and Mazur over nine decades ago.
For over four decades this definition is considered also for a arbitrary convex bodies of $E^d$.
Namely, for convex bodies $C, D$ of $E^n$ this extended Banach-Mazur distance sounds as follows

$$\delta_{\rm BM}(C, D) = \inf_{a,\, h_{\lambda}} \{ \lambda;\ \ a(D) \subset C \subset h_{\lambda}a(D) \big).$$

\noindent
where $a$ stands for an affine transformation and $h_\lambda$ denotes a homothety with a positive ratio $\lambda$.
For the relationship of them see Claim of \cite{[L5]}.

A survey on the Banach-Mazur distance is given in the book \cite{[TJ]} by Tomczak-Jaegerman.
Moreover, in Sections 3.2 and 3.3 of the book \cite{[T]} by Toth, and in Section 4.1 of the book \cite{[AS]} by Aubrun and Szarek. 

Here is the notion of the centroid Banach-Mazur distance of convex bodies $C, D$ of $E^d$:

$$\delta_{\rm BM}^{\rm cen} (C,D) = \inf_{a,\, h_\lambda} \{\lambda ; \, a(D) \subset C \subset h_\lambda a(D) \ {\rm and} \ {\rm cen}(a(D)) =  {\rm cen}(C) \},$$

\noindent
where $a$ again stands for an affine transformation, but $h_\lambda$ means a homothety with the ratio $\lambda \geq 1$ whose center is at the centroids ${\rm cen}(a(D)) =  {\rm cen}(C)$ of $a(D)$ and $C$. 
Observe that the centroids of $C, D$ in this definition take over the roles of the centers of the centrally-symmetric bodies in the original definition of Banach-Mazur distance.
We easily show that $\delta_{\rm BM}^{\rm cen} (C,D) = \delta_{\rm BM}^{\rm cen} (D,C)$ for every $C$ and $D$. 
Recall that pioneer research on the centroid was provided by Neumann \cite{[N]}.

In Theorem we prove that $\delta_{\rm BM}^{\rm cen} (P,T) = \frac{5}{2}$ for the parallelogram $P$ and the triangle $T$ in $E^2$.
Our effort is put in order to show that $\delta_{\rm BM}^{\rm cen} (P,T) \geq \frac{5}{2}$ since the opposite inequality immediately follows from easy examples.

At the end of the paper we present a few remarks.
The first concerns the positions of our triangle with respect to the parallelogram for which the ratio $\frac{5}{2}$ is realized. 
The second comments the dual version of Theorem.
The third shows the following generalization of Theorem for a centrally-symmetric convex body $M$ in place of $P$: for every triangle $T$ inscribed in $M$ with the common centroid we have $M \subset 3T$.
We also propose a more general task to consider an arbitrary convex body instead of $M$.
Finally, we ask about a generalization of Theorem for $E^d$.

As usual, by $\inter (A)$ and $\bd (A)$ we denote the interior and boundary of a set $A$.

\section{The distance between the parallelogram and the triangle is $\frac{5}{2}$}

\vskip0.1cm
\noindent
{\bf Theorem.}
{\it For the parallelogram $P$ and the triangle $T$ we have $\delta_{BM}^{\rm cen} (P, T) = \frac{5}{2}$.} 

\vskip-0.8cm
{\ }

\begin{proof}
In the proof, by $\lambda C$ we mean the homothetic image of a set $C$ with a positive ratio $\lambda$ and the center at the origin $o$ of $E^2$.
By $S$ denote the square with vertices $(1 , \pm 1)$ and $(-1 , \pm 1)$.

Let us rephrase our theorem as the conjunction of the two following sentences 

(*) for every triangle $\Delta \subset S$ with the centroid at the center $o$ of $S$ the interior of the triangle $\frac{5}{2}\Delta$ does not contain $S$, 

(**) there exists a triangle $\Delta_0 \subset S$ with the centroid in the center of $S$ for which $\frac{5}{2}\Delta_0$ contains $S$.  

Taking as $\Delta_0$ the triangle with vertices $(1, \frac{1}{2}), (-1, \frac{1}{2})$ and $(0, -1)$ we see that (**) is true.

Later, our aim is to show that (*) holds true.

\vskip0.1cm
The following obvious fact is applied soon two times in the proof.

\vskip0.1cm
(\hskip -0.7cm \begin{tabular}{ll}
\rotatebox{45}{} &\rotatebox{45}{$\Delta$}  \\
\end{tabular} \hskip -0.35cm){\it \hskip0.2cm If} (*) {\it is true for a triangle centered at $o$ and containing $\Delta$, then} (*) {\it is true for $\Delta$.}

\vskip 0.2cm

By 
(\hskip -0.7cm \begin{tabular}{ll}
\rotatebox{45}{} &\rotatebox{45}{$\Delta$}  \\
\end{tabular} \hskip -0.35cm){\it \hskip0.02cm} 
we may assume that $\Delta \subset S$ has at least one vertex in the boundary of $S$.
Still if all the vertices of the triangle are in the interior of $\inter (S)$, we can increase this triangle by a homothety with center $o$ such that at least one of its vertices ``arrives" to the boundary of $S$.

Without losing the generality assume that such a vertex $a$ is $(1,\alpha)$, where $0 \leq \alpha \leq 1$.

Denote by $b$ and $c$ the endpoints of the opposite side of $\Delta$ such that $a, b, c$ are in the positive order.
The midpoint of $bc$ is $d=(-\frac{1}{2}, -\frac{1}{2}\alpha)$.
The reason is that the segment connecting $d$ with $a$ passes through $o$ which is the centroid of $\Delta$.  

We may assume that $b$ or $c$ is in $\bd (S)$.
Still if $b$ and $c$ are in $\inter (S)$, we do not loose the generality by properly enlarging $\Delta$.
Namely, we increase $bc$ by the homothety at its center $d$ as long as an endpoint attains the boundary of $S$, but clearly, the other must remain in $S$.
This follows by 
(\hskip -0.7cm \begin{tabular}{ll}
\rotatebox{45}{} &\rotatebox{45}{$\Delta$}  \\
\end{tabular} \hskip -0.35cm){\it \hskip0.02cm} 
and the observation that the center of the increased triangle is $o$ again.

\vskip 0.3cm
{\bf Case 1}, when the vertex $b$ attains the boundary of $S$ not later than $c$.

{\bf Subcase 1.1}, when $\alpha \in (0,1]$. 

Since $\alpha >0$, we have $- \frac{1}{2} \alpha <0$.
Thus $d$ is below the axis $y=0$. 
Since $d$ is the midpoint of $bc$, this and the fact that $c$ is over or on the straight line $y=-1$ imply that $b$ must be below the straight line $y=1$.
Hence $b$ does not belong to the side $(1,1)(-1,1)$.
Obviously, it also does not belong to the sides $(-1,-1)(1,-1)$ and $(1,-1)(1,1)$.
So $b$ belongs to the side $(-1,1)(-1,-1)$ and has the form $b = (-1, \beta)$.
See Figure 1.
Clearly, $-1 \leq \beta \leq 1$.

Applying again the fact that $d=(-\frac{1}{2}, -\frac{1}{2}\alpha)$ is the middle of $bc$, we conclude that $c= (0, -\alpha -\beta)$.
Since the second coordinate of $c$ is at least $-1$, we get $\alpha-\beta \geq -1$ which means that $1-\alpha \geq \beta$. 
Since the order of vertices $a, b, c$ of $\Delta$ is positive, we get $\beta \geq -\alpha$ (for $\beta = -\alpha$ our triangle $\Delta$ degenerates to a segment).
Resuming, $-\alpha \leq \beta \leq 1-\alpha$ and thus $b \in gh$, where $g = (-1, 1-\alpha)$ and $h = (-1,-\alpha)$.
 
Take into consideration the triangle $\frac{5}{2}\Delta$.
Its corresponding vertices are $a'= (\frac{5}{2}, \frac{5}{2}\alpha)$, $b' = (-\frac{5}{2}, \frac{5}{2}\beta)$ and $c' =(0, -\frac{5}{2 }\alpha - \frac{5}{2}\beta)$.
Here are the equations of the straight lines $\ell_{a'c'}$ and $\ell_{b'c'}$ containing the sides $a'c'$ and $b'c'$, respectively:

$\ell_{a'c'}$: $y- \frac{5}{2}\alpha = (2\alpha + \beta)(x- \frac{5}{2})$,

$\ell_{b'c'}$: $y+ \frac{5}{2}\alpha +  \frac{5}{2} \beta = (-\alpha - 2\beta)x$. 

\vskip0.1cm
We intend to show that for every $\alpha \in (0,1]$ and $\beta \in [-\alpha, 1-\alpha]$ at least one of the straight lines $\ell_{a'c'}$, $\ell_{b'c'}$ intersects the square $S$. 
This would mean that the interior of $\frac{5}{2}\Delta$ does not contain the whole $S$.
Consequently (*) would hold true.

Consider the parallelogram 
$V= \{(\alpha, \beta); \ 0< \alpha \leq 1, \alpha \geq \beta \geq 1 - \alpha\}$
in the coordinate system $o\alpha\beta$, where $o =(0,0)$, without the side $(0,0)(0,1)$. 
See Figure 2.

The below consideration shows that for every $(\alpha, \beta) \in R$ at least one of the following two sentences is true in the system $Oxy$:

$\ell_{a'c'}$ intersects the side $(-1,-1)(1,-1)$, 

$\ell_{b'c'}$ intersects the side $(-1,-1)(-1,1)$.

\vskip 0.1cm
We intersect $\ell_{a'c'}$ with the line $y=-1$ containing the side $(-1,-1)(1,-1)$.
The point of intersection is $e = \Big(\frac{2 - 5\alpha + 5\beta}{4\alpha + 2\beta}, -1 \Big)$.
If the line $\ell_{a'c'}$ intersects the side $(-1,-1)(1,-1)$, from $x \leq 1$ we obtain $\beta \leq \frac{2-\alpha}{3}$.
As a result we see the region $V(\ell_{a'c'})$ of points $(\alpha, \beta)$ in the coordinate system $o\alpha\beta$ for which the line $\ell_{a'c'}$ intersects the side $(-1,-1)(1,-1)$; it is the part of $R$ not above the straight line $\beta = \frac{2-\alpha}{3}$. 

The straight line $\ell_{b'c'}$ intersects the one $x=-1$ containing the side $(-1,1)(-1,-1)$. 
The point $f$ of intersection is $(-1, -\frac{3}{2}\alpha -\frac{1}{2}\beta)$.
If the line $\ell_{a'c'}$ intersects the side $(-1,1)(-1,-1)$, then from $y \geq -1$ we conclude that $\beta \leq 2 -3\alpha$.
The region $V(\ell_{b'c'})$ of points $(\alpha, \beta)$ in the system $o\alpha\beta$ for which the line $\ell_{b'c'}$ intersects the side $(-1,-1)(-1,1)$ is the part of $R$ not above the straight line $\beta = 2-3\alpha$.

\eject

\vskip0.3cm
\begin{figure}[htbp]        
\hskip0.2cm
\includegraphics[width=10.35cm]{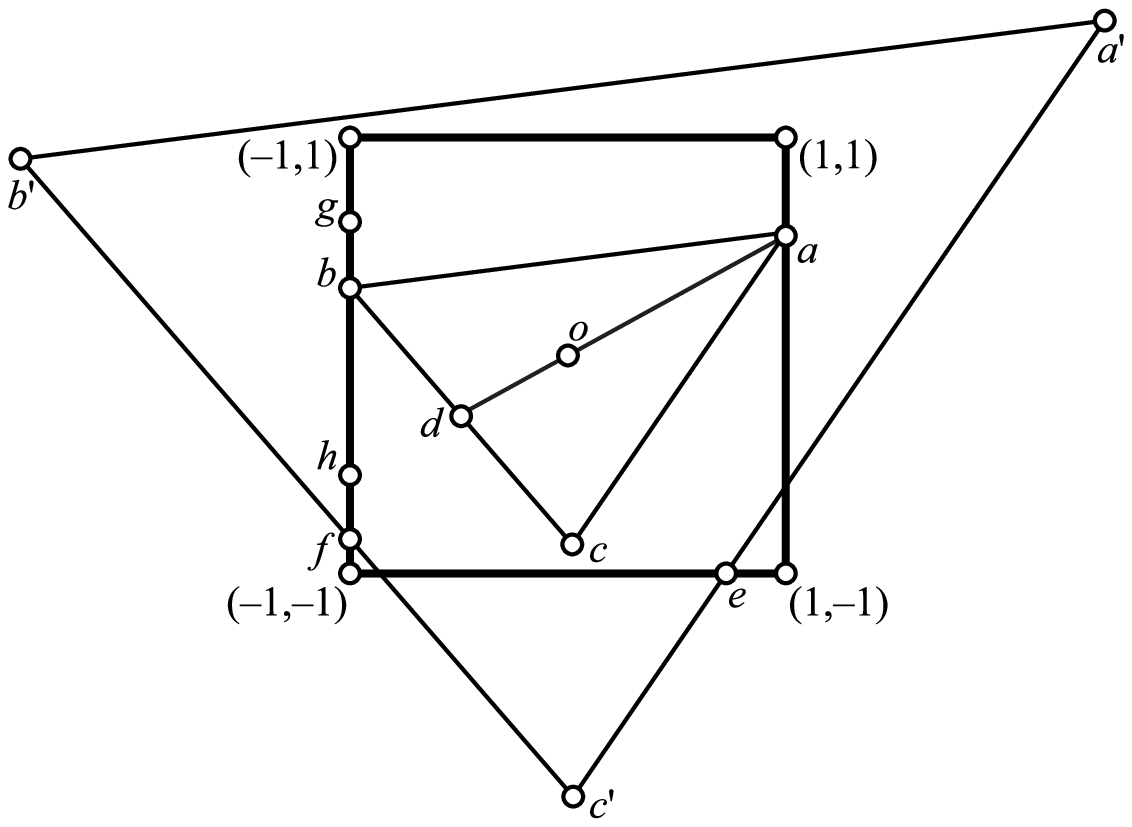} \hskip-0.1cm 
\includegraphics[width=4.95cm]{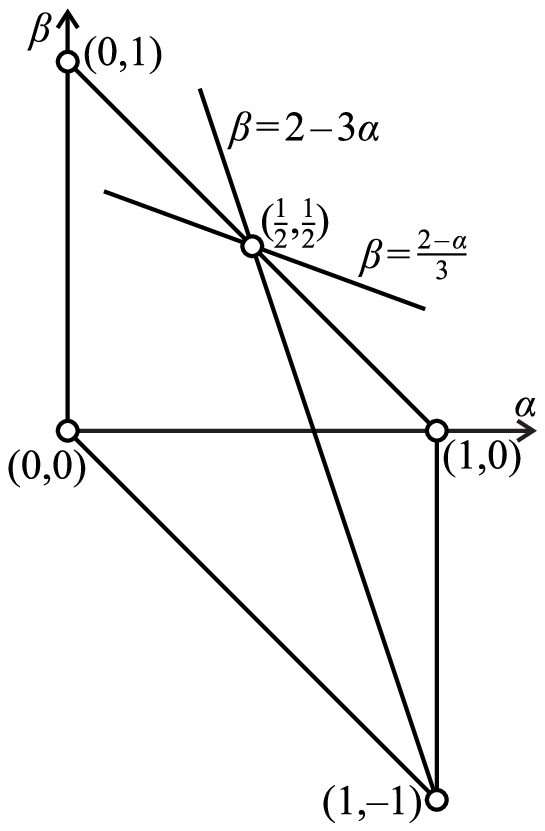}  \\ 

\vskip-0.2cm
\hskip 2.5cm {\bf Fig. 1.} Positions of $a, b, c$ in Case 1 \hskip 2.7cm 
{\bf Fig. 2.} Region $V$ 

\end{figure}

\vskip-0.3cm
These two straight lines $\beta = \frac{2-\alpha}{3}$ and $\beta = 2-3\alpha$ intersect at the point $(\frac{1}{2}, \frac{1}{2})$ of the system $o\alpha\beta$ which is in the side $(1,0)(0,1)$ of $R$.
From this and the two preceding paragraphs we see that $V \subset \ell_{a'c'} \cup V(\ell_{b'c'})$. 
Hence (*) holds true.

\vskip 0.2cm
{\bf Subcase 1.2}, when $\alpha=0$.

We have $d=(-\frac{1}{2}, 0)$ which implies that $b \in (0,1)(-1,1)$ and $c \in (-1,-1)(0,-1)$ are symmetric with respect to $d$.
We let the reader to check that
the interior of the triangle $\frac{5}{2}\Delta$ does not contain $S$, so (*) holds true.

\vskip 0.3cm
{\bf Case 2}, when the vertex $c$ attains the boundary of $S$ not later than $b$.

\vskip0.2cm
{\bf Subcase 2.1} when $\alpha \in (0, 1]$.

Observe that $c$ must be in the side $(-1,-1)(1,-1)$ (see Figure 3).
So $c$ has the form $(\gamma, -1)$, where $-1 \leq \gamma \leq 1$.
From $b \in S$ and the fact that $d$ has the first coordinate $-\frac{1}{2}$ we see that $-1 \leq \gamma \leq 0$.
Since $d$ is the midpoint of $bc$, we have $b = (-1-\gamma, 1-\alpha)$.

Take into account the triangle $\frac{5}{2}\Delta$.
Its corresponding vertices are 
$a'= (\frac{5}{2}, \frac{5}{2}\alpha)$, 
$b' = (-\frac{5}{2}- \frac{5}{2}\gamma, \frac{5}{2}- \frac{5}{2}\alpha)$ and 
$c' =(\frac{5}{2}\gamma, -\frac{5}{2})$ (see Figure 3).

Here are the equations of the straight lines containing the sides of the triangle $a'b'c'$. 

$\ell_{a'b'}:  y - \frac{5}{2}\alpha = \frac{-1+2\alpha}{2+\gamma}(x - \frac{5}{2})$, 

$\ell_{b'c'}:  y + \frac{5}{2} = \frac{2-\alpha}{-1-2\gamma}(x - \frac{5}{2}\gamma)$, 

$\ell_{a'c'}:   y - \frac{5}{2}\alpha = \frac{\alpha+1}{1-\gamma}(x - \frac{5}{2})$.

We intend to show that for every $\gamma \in [-1, 0]$ at least one of these lines $\ell_{a'b'}$,  $\ell_{b'c'}$ and $\ell_{a'c'}$ intersects the square $S$. 
This would mean that (*) holds true.

\eject

\vskip0.3cm
\begin{figure}[htbp]        
\hskip0.1cm
\includegraphics[width=10.12cm]{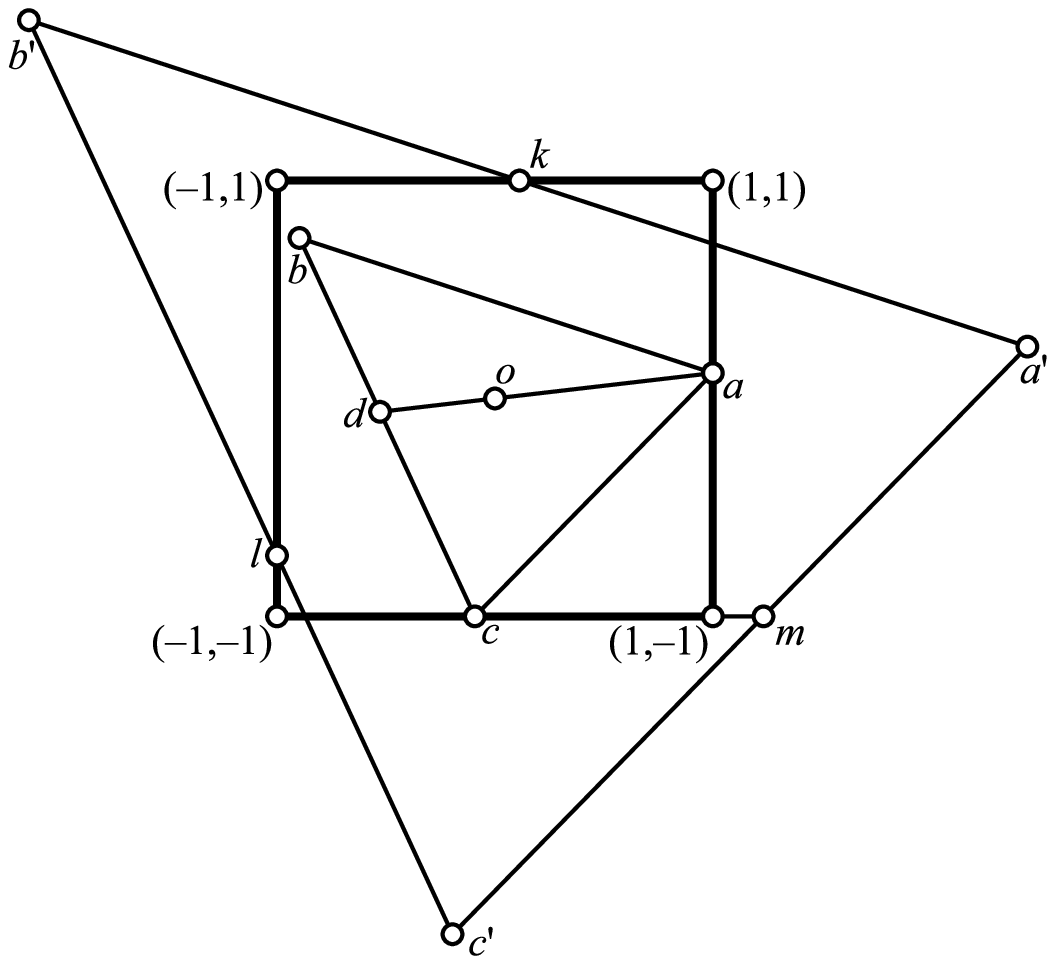} \hskip-0.1cm 
\includegraphics[width=5.984cm]{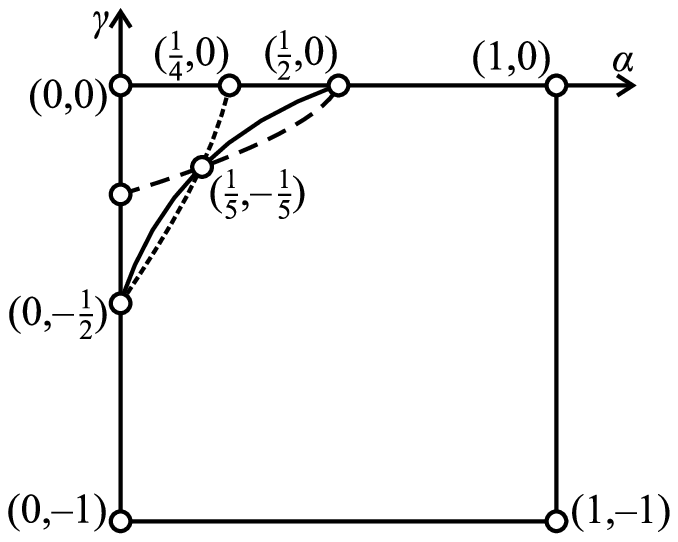}  \\ 

\vskip-0.3cm
\hskip 1.7cm {\bf Fig. 3.} Positions of $a, b, c$ in Case 2 \hskip 3.85cm 
{\bf Fig. 4.} Region $W$ 

\end{figure}

\vskip-0.1cm
Now let us deal with the set $W= \{(\alpha, \gamma): 0 < \alpha \leq 1, -1 \leq \gamma \leq 0 \}$ of points in the coordinate system $o\alpha\gamma$ (see Figure 4), where $o = (0,0)$. 

The point $k$ of intersection of $\ell_{a'b'}$ with the straight line $x=1$ containing the side  $(1,1)(1,-1)$ has the second coordinate $y= -\frac{3}{2} \cdot \frac{-1+2\alpha}{2+\gamma} + \frac{5}{2}\alpha$. 
If $\ell_{a'b'}$ intersects this side, then from $y \leq 1$ we get $-\frac{3}{2} \cdot \frac{-1+2\alpha}{2+\gamma} + \frac{5}{2}\alpha \leq 1$.
Equivalently, $\gamma \geq \frac{1- 4\alpha}{-2+5\alpha}$.
So the curve $\gamma = \frac{1- 4\alpha}{-2+5\alpha}$ in the coordinate system $o\alpha\gamma$ intersects the axis $o\alpha$ at $(\frac{1}{4}, 0)$ and the axis $o\gamma$ at $(0, -\frac{1}{2})$.
This permits to see the subregion $W(\ell_{a'b'})$ of $W$ of points for which $\ell_{a'b'}$ intersects $S$. 
It is bounded by the piece of $\gamma = \frac{1- 4\alpha}{-2+5\alpha}$ for $0\leq \alpha \leq \frac{1}{4}$ (marked by the dotted line in Figure 4) and the segments $(0,0)(0, \frac{1}{4})$ and $(0,0)(0,-\frac{1}{2})$ in the system $o\alpha\gamma$.

The intersection of $\ell_{b'c'}$ with the straight line $y= -1$ containing the side $(-1,-1)(1,-1)$ is at a point $l=(x, -1)$, 
where $x$ fulfills $\frac{3}{2} = \frac{2-\alpha}{-1-2\gamma}(x - \frac{5}{2}\gamma)$, this is
$x= \frac{-3 -6\gamma}{4-2\alpha} +\frac{5}{2}\gamma$.
If $\ell_{b'c'}$ intersects the side $(1,-1)(-1,-1)$, then from $x \geq -1$ we obtain $-1 \leq \frac{-3 -6\gamma}{4-2\alpha} +\frac{5}{2}\gamma$ which (for $\alpha \in [0, \frac{1}{2}]$) is equivalent to $\gamma \geq \frac{1-2\alpha}{-4+5\alpha}$.
We easily check that the curve $\gamma = \frac{1-2\alpha}{-4+5\alpha}$ in the coordinate system $o\alpha\gamma$ intersects the axis $o\alpha$ at $(\frac{1}{2},0)$ and the axis $o\gamma$ at $(0, -\frac{1}{4})$.
In Figure 4 we see the subregion $W(\ell_{b'c'})$. 
It is bounded by the piece of $\gamma = \frac{1-2\alpha}{-4+ 5\alpha}$ for $0\leq \alpha \leq \frac{1}{2}$ (marked by the dashed line in Figure 4)  and the segments $(0,0)(0, \frac{1}{2})$ and $(0,0)(0,-\frac{1}{4})$ in the coordinate system $o\alpha\gamma$.

The intersection of $\ell_{a'c'}$ with the straight line $y= -1$ is at a point $m= (x, -1)$, where $x$ fulfills
$-1 - \frac{5}{2}\alpha = \frac{\alpha+1}{1-\gamma}(x - \frac{5}{2})$,
this is $x= \frac{5}{2} + \frac{-2 -5\alpha}{2+2\alpha} (1- \gamma)$.
If $\ell_{a'c'}$ intersects the side $(-1,-1)(1,-1)$, than from $x \leq 1$ we obtain $\frac{3}{2} \leq \frac{2+5\alpha}{2+\alpha} (1- \gamma)$ which, for our positive $\alpha$ is equivalent to $\gamma \leq \frac{-1+2\alpha}{2+5\alpha}$.
We easily check that the curve $\gamma = \frac{-1+2\alpha}{2+5\alpha}$ in the coordinate system $o\alpha\gamma$ intersects the axis $o\alpha$ at $(\frac{1}{2},0)$ and the axis $o\gamma$ at $(0, -\frac{1}{2})$.
In Figure 4 we see this subregion $W(\ell_{a'c'})$ of $W$ of points for which $\ell_{a'c'}$ intersects $S$.  
This subregion is bounded by the piece of the curve $\gamma = \frac{-1 +2\alpha}{2+ 5\alpha}$ for $0 < \alpha \leq \frac{1}{2}$ (marked by the solid line in Figure 4) and by the segments connecting the succeeding pairs of points $(0, -\frac{1}{2}), (0,-1), (1,-1), (1,0)$ and $(\frac{1}{2}, 0)$ in the coordinate system $o\alpha\gamma$.
Clearly, here the first of these segments is not in $W(\ell_{a'c'})$.

Observe that the point $(\frac{1}{5}, -\frac{1}{5})$ belongs to the three pieces of curves bounding our three considered subregions.
Moreover, $\frac{-1+2\alpha}{2+5\alpha} \geq \frac{1-4\alpha}{-2+5\alpha}$ for $0< \alpha \leq \frac{1}{5}$ 
and $\frac{-1+\alpha}{2+5\alpha} \geq \frac{1-2\alpha}{-4+5\alpha}$ for $\frac{1}{5} \leq \alpha \leq \frac{1}{2}$ (see Figure 4).
These facts and the three preceding paragraphs imply that $W \subset W(\ell_{a'b'}) \cup W(\ell_{b'c'}) \cup W(\ell_{a'c'})$.

\vskip0.2cm
{\bf Subcase 2.2} for $\alpha=0$.  
 
Now we have $d=(-\frac{1}{2}, 0)$.
This implies that $c \in (-1,-1)(0,-1)$ and $b \in (0,1)(-1,1)$ are symmetric with respect to $d$.
It is easy to check that the interior of the triangle $\frac{5}{2}\Delta$ does not contain $S$.
So (*) holds true.

\vskip 0.2cm
From Cases 1 and 2 we see that always  $\inter(\frac{5}{2}\Delta)$ does not contain $S$, i.e., (*) is fulfilled.
\end{proof}

The fact proved in Theorem is claimed without a proof at the bottom of p. 259 of \cite{[G]}.

\section{Final remarks}

Let us comment the positions of the triangle $\Delta_0$ in $S$ (drawn by a thick line in Figure 5) from the proof of Theorem, for which (**) holds true.
As it follows from our considerations, the only two such triangles $\Delta_0 \subset S$ (up to symmetric positions) are the following.
The first has vertices $a_1 = (1, \frac{1}{2})$, $b_1 = (-1, \frac{1}{2})$ and $c_1 = (0, -1)$ (by the way, two symmetric positions are also seen for $\alpha =0$ in Subcases 1.2 and 2.2 of the proof of Theorem). 
The second has vertices $a_2 = (1, \frac{1}{5})$, $b_2 = (-\frac{4}{5}, \frac{4}{5})$ and $c_2 = (-\frac{1}{5}, -1)$. 
We see both and also corresponding $\frac{5}{2}\Delta_0$ in Figure 5.

Look at the dual situation.
Putting $C=P$ and $D=T$ in $\delta_{\rm BM}^{\rm cen} (C,D) = \delta_{\rm BM}^{\rm cen} (D,C)$, by Theorem we get $\delta_{\rm BM}^{\rm cen} (T,P) = \frac{5}{2}$. 
In particular, for a given equilateral triangle $T$ (marked by a thick line in Figure 6) there are only two extreme positions (up to symmetric ones) of a parallelogram $P_0$ for which $\frac{5}{2}P_0$ contains $T$.

In connection to our Theorem recall Theorem from \cite{[L3]} shows that every centrally symmetric convex body $M \subset E^2$ permits to inscribe a triangle $\Delta$ whose centroid is at the center of symmetry of $M$ such that $M \subset \frac{5}{2}\Delta$ (our present Theorem explains that this ratio cannot be lessened when $M$ is a parallelogram.)
Here is a claim which considers any inscribed triangle.

\vskip0.3cm
{\bf Claim.}
{\it Let $M \subset E^2$ be a centrally-symmetric convex body. 
For every inscribed triangle $\Delta$ in $M$ with the common centroid we have $M \subset 3\Delta$.} 

\vskip-0.9cm
\begin{proof} 
Let $\Delta= abc$. 
Take the symmetric triangle $\Delta_s = a_sb_sc_s$ with respect to $o$.
The convex

\eject

\vskip0.2cm
\begin{figure}[htbp]        
\includegraphics[width=9.8cm]{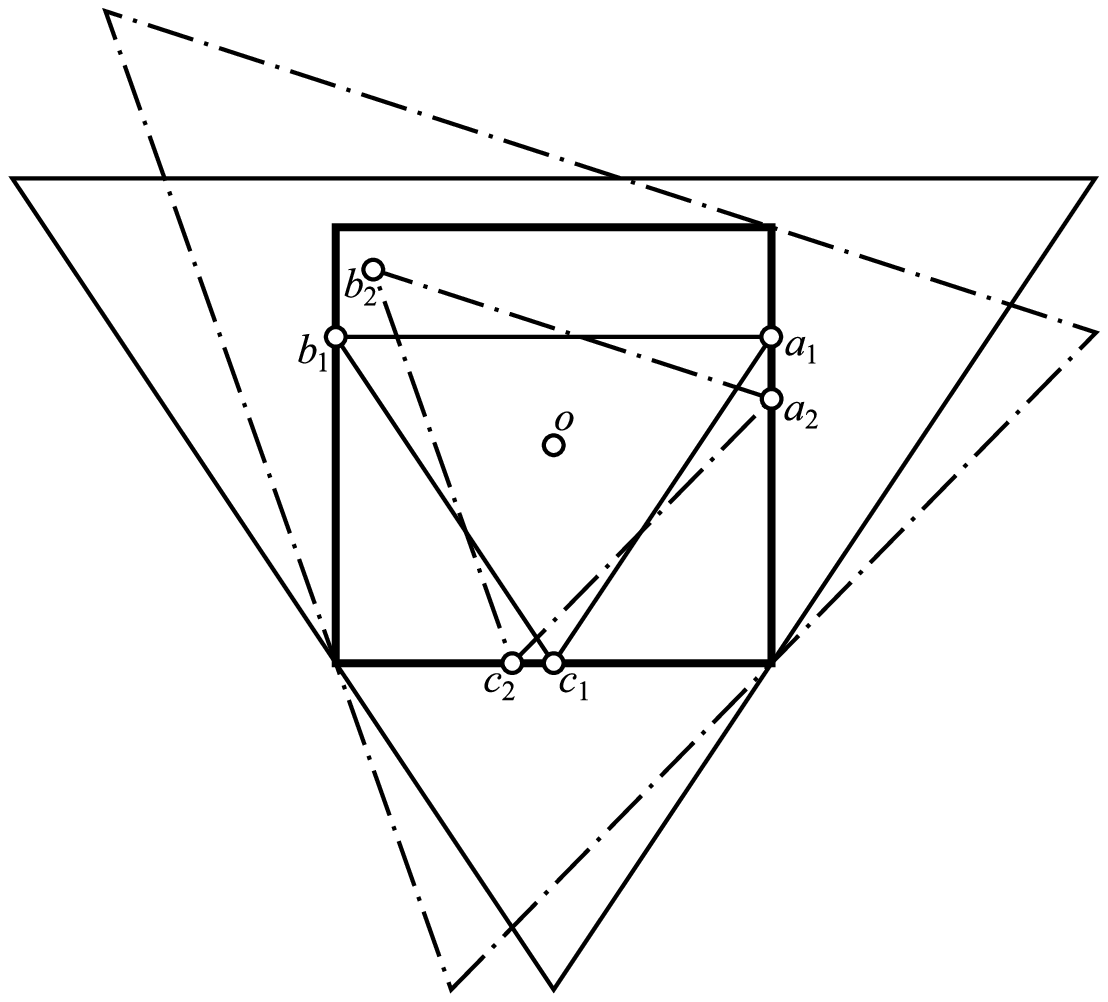} \hskip-0.3cm 
\includegraphics[width=6.1cm]{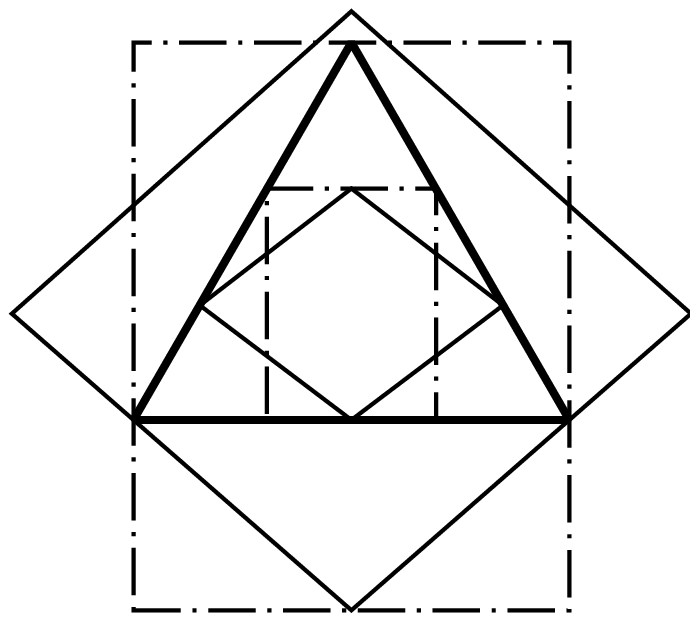}  \\ 

\hskip 1.4cm {\bf Fig. 5.} Extreme positions of $\Delta$ in $S$ \hskip 1.75cm 
{\bf Fig. 6.} Extreme positions of $P$ in $T$

\end{figure}

\noindent
hull $H$ of $\Delta \cup \Delta_s$ is an affine-regular hexagon inscribed in $M$.
Prolonging the sides $ab_s, bc_s, ca_s$ we get three points of intersection.
Denote by $S(H)$ the star being union of the triangle with these three vertices and the symmetric triangle with respect to $o$.
By the convexity of $M$ we conclude that $M \subset S(H)$.
Moreover, observe that $S(H) \subset 3\Delta$.
Consequently, $M \subset 3 \Delta$.
\end{proof}

By the way, Figure 3 of \cite {[L1]} shows an analogous situation of extreme homothetic parallelograms $S \subset \Delta$ and $(\sqrt 2 +1)S \supset \Delta$ without the requirement that centroids of $S$ and $\Delta$ coincide.
The ratio $3$ in Claim cannot be lessened as it shows the example of the square with vertices $(-1,\pm 1), (1,\pm 1)$ and the triangle with vertices $(1,1), (-1,0), (0,-1)$. 

Theorem of  \cite{[L3]} says the following. 
{\it Let $M \subset E^2$ be a centrally symmetric convex body.
In $M$ it can be inscribed a triangle $\Delta$ whose centroid is the center of symmetry of $M$ such that $M \subset \frac{5}{2}\Delta$.}
From our Theorem it follows that this ratio $\frac{5}{2}$ cannot be lessened.

What about considering an arbitrary convex body in place of $M$ in Claim?
The author conjectures that for every planar convex body $C$ and an arbitrary inscribed triangle $\Delta$ with the common centroid we have $C \subset 4\Delta$.
Possibly the result of Neumann \cite{[N]} would be a good tool.
The the ratio $4$ cannot be lessened for $C$ being a triangle with the centroid at $o$ and $\Delta = -\frac{1}{4}C$.

The above remarks can be seen in the wider context of approximation by triangles in \cite{[BL]}.

Finally, let us ask about a higher dimensional generalization of Theorem.
In $E^3$ the paralellotope $Q$ permits to inscribe (thus also put inside) a simplex $S$ whose centroid is in the center of $Q$ and $Q \subset 3S$.
Just put the vertices of $S$ at some four non-neighboring vertices of $Q$. 
Thus $\delta_{BM}^{\rm cen} (Q, P) \leq 3$.
The author believes that the equality holds true here, but possible evaluation seems to be complicated.
For higher dimensions the task of finding or at least estimating $\delta_{BM}^{\rm cen} (Q, P)$ remains open. 
From \cite{[L4]} we only conclude that $\delta_{BM}^{\rm cen} (Q, P) \leq 2n-1$.
Since $\delta_{BM}^{\rm cen}(P, Q) = \delta_{BM}^{\rm cen} (Q, P)$, the same also follows from 
the last paragraph of \cite{[L2]}.

\end{document}